\documentclass[11pt, oneside]{article}   	
\usepackage{amsmath}
\usepackage{mathtools} 				
\usepackage{geometry}                		
\usepackage{graphicx}				
\usepackage{amsmath}								
\usepackage{amsthm}
\usepackage{amssymb}
\usepackage{tikz}
\usepackage{float}
\usepackage{enumitem}
 \usepackage{tabu}
 \usetikzlibrary{decorations.pathreplacing}
\newtheorem{theorem}{Theorem}
\newtheorem*{main theorem}{Main Theorem}
\newtheorem{corollary}{Corollary}

\newtheorem{lemma}{Lemma}
\theoremstyle{definition}	
\newtheorem{defn}{Definition}	
\newtheorem{empl}{Example}
\newtheorem{rmk}{Remark}

\newcommand{\gen}{\text{gen}}
\usepackage[utf8]{inputenc}

\DeclarePairedDelimiter{\floor}{\Big \lfloor}{\Big \rfloor}

\title{The generic value set of a two-generator semigroup}
\author{Justin Lake \& Lee McEwan}

\begin{document}
\maketitle

\begin{abstract}
This article gives a recursion for the minimal generators $G=\{g_i\}$ of the generic value set $\Lambda_{\gen}$ of a plane curve germ $C$ with a two-generator semigroup $\Gamma \langle p, m \rangle$. The main result provides for explicit calculation of all $g_i$ and shows that they are in fact minimal generators. Explicit formulas for the cardinality of the minimal generators $|G|$ and for the conductor $c(\Lambda_{\gen})$ of $\Lambda_{\gen}$ are given.  The recursion can be used to compute the generic Tjurina number $\tau_{\gen}$, a method compared to that given in \cite{ABM}. The main result is based on the non-explicit algorithm and ideas of Delorme \cite{D1}, \cite{D2}.
\end{abstract}

\section{Introduction} \label{intro}
Given a singular plane branch, i.e., the germ of an irreducible singular plane curve, one can ask: what are the valuations of the set of differential one-forms? Fixing a topological class of singular plane branches, there
is a generic set of such valuations, those appearing for a generic branch of the specified class.
In this paper we give explicit formulas for the generators of this generic set, in the case where the 
Puiseux characteristic is $(p;m)$, with $1<p<m$ and $\operatorname{gcd}(p,m)=1$.

Zariski \cite{Z} shows that any two singular plane branches of the same topological type will have the same Puiseux characteristic, equivalently the same associated  numerical semi-group $\Gamma$. If the branch is given by $f(X,Y)=0$, with $\mathcal{O}:=\mathbb{C}[[X,Y]]/f$, then $\Gamma=v(\mathcal{O})$,
where $v$ is the valuation for functions.
For this fixed topological type there may be different curves that have different sets of valuations of differential one-forms, which we will call \emph{value sets}; a value set will be denoted by $\Lambda$. 
Thus $\Lambda=v(\mathcal{O}d\mathcal{O})$, where $v$ now denotes valuation for one-forms.
In defining $v$, we adopt the convention that if $g$ is a function vanishing at the singular point, then
$v(dg)=v(g)$. Thus, for example, for the branch with $X=t^2$ and $Y=t^3$, we have $v(YdX)=v(Y)+v(X)=5$;
equivalently $v(YdX)=v(2t^4dt)=v(2t^4)+1$. This is in accord with the convention of Hefez and Hernandes \cite{HH}.
With this convention, the value set $\Lambda$ contains $\Gamma^*:=\Gamma\setminus\{0\}$.
Also note that for any $r\in\Lambda,$ we have $\Gamma+r\subset\Lambda,$ making $\Lambda$ a $\Gamma$-semimodule. 

Fixing the topological type, there is a generic set of curves that have the same value set, denoted $\Lambda_{\gen}.$ Peraire \cite{P} provides a difficult algorithm that gives the generic value set for a general Puiseux characteristic. In a brief article, Delorme \cite{D1} gives a much simpler algorithm for computing the generators of $\Lambda_{\gen}$ for the case where the Puiseux characteristic is $(p;m)$, equivalently
$\Gamma=\langle p,m\rangle.$ While the statement of this algorithm is simple, computing with it is impractical; an example is provided to illustrate the challenge. In this paper we make use of Delorme's algorithm to compute formulas for the minimal generators for $\Lambda_{\gen}$ as a $\Gamma$-semimodule in the case where $\Gamma=\langle p,m\rangle.$ Our formulas are easy to execute by hand even for large values of $p$ and $m$. From our basic result we also deduce simple formulas for the conductor of $\Lambda_{\gen}$ and the cardinality of the set of minimal generators. Lastly, we show how the minimal Tjurina number for given $\langle p,m\rangle$ can be easily calculated from our result. This calculation appears to be quite different from the (equivalently easy to use, and more general) formula in \cite{ABM}. 

We begin in Section~\ref{DelormeAlgorithm} by stating Delorme's algorithm.
Sections~\ref{prelims} and \ref{MT} present our recursive formulas
for the generators produced by the algorithm: first we lay out some preliminaries,
and then we state and prove the main theorem (Theorem  \ref{Main}).
These formulas typically compute precisely a minimal set of 
generators for the generic value set. The simple situation in which a few non-minimal generators may be produced is described in Section~\ref{specialcases}.
Section~\ref{examples} is devoted to examples, and
Section~\ref{Tjurinaformula} relates our work to
the formula for the Tjurina number found in \cite{ABM}.

\section{Delorme's algorithm}\label{DelormeAlgorithm}

Here we present the algorithm from Delorme's 1974 article \cite{D1}. This will be the process we use to prove that we are indeed working in the generic case, by taking what Delorme calls minimal steps. We must first define some notation.

\begin{defn} \label{Eandu}
Let $\Gamma$ be a numerical semigroup, and let $g_{-1}<g_0<g_1<\dots<g_n$ be the minimal set of generators of a $\Gamma$-semimodule $\Lambda$. 
For $-1\leq i\leq n$, define
\begin{equation}
\label{Ei}
E_i:=\bigcup_{-1\leq j\leq i}(\Gamma+g_j).
\end{equation}
Because $\{g_i\}$ is the {\em minimal} generating set, we have $g_i\notin E_{i-1}$. Clearly $E_n=\Lambda$.
We also set  (for $0\leq i\leq n$)
\begin{equation}
\label{defui}
u_i:=\operatorname{min}\left((\Gamma+g_i)\cap E_{i-1}\right). 
\end{equation}
\end{defn}

\begin{defn} 
Let $\Gamma$ be a numerical semigroup. The {\em conductor} of $\Gamma$ is
\begin{equation}
\label{Cond}
\mu:=\operatorname{min}\{x \in \Gamma \mid \mathbb{N}+x \subset \Gamma\}.
\end{equation}
Similarly, let $\Lambda$ be a value set for $\Gamma$. The {\em conductor} of $\Lambda $ is
\begin{equation}
\label{CondV}
c(\Lambda):=\operatorname{min}\{x \in \Lambda \mid \mathbb{N}+x \subset \Lambda\}.
\end{equation}
\end{defn}

We remark that $\mu$ is also called the {\em Milnor number} (\cite{M}, \cite{HR}). For the remainder of the paper, we assume $\Gamma$ is of the form $\Gamma=\langle p,m\rangle$. In this case it is well-known that the conductor of $\Gamma$ is given by $\mu = (p-1)(m-1)$. 

We now present Delorme's algorithm for obtaining the generic $\Gamma$-semimodule $\Lambda_{\gen}.$
\\
\\
{\bf Delorme's algorithm:}
Begin with $g_{-1}=p$ and $g_0=m$;
this gives $E_{-1}=\Gamma+p$ and
$E_0=\Gamma^*$.
For $i \geq 1$ define 
\begin{equation}
\label{gi}
g_i:=\text{min}\left((\mathbb{N}+u_{i-1})\setminus E_{i-1}\right).
\end{equation}
Using formulas (\ref{Ei}) through (\ref{gi}), for each index $i$ compute first the value of $g_i$, then $E_i$, then $u_i$. Because $u_i\in E_{i-1}$ and $g_{i+1} \not\in E_i$, the process ensures $g_i < u_i < g_{i+1}$. 
Because $\{u_i\}$ is strictly increasing and $\Gamma$ has finite complement in $\mathbb{N}$, there exists a smallest positive $i = n$ such that $\mathbb{N} + u_n \subset E_n$. Then the set (\ref{gi}) is empty and the algorithm stops. The values $g_{-1}, g_0, g_1, \dots, g_n$ computed up to this stopping point are the minimal generators of $\Lambda_{\gen}.$

\begin{rmk} The conductor of $\Gamma$ provides a crude upper bound on the number of steps in the algorithm.
Ultimately we will find a precise value for the number of minimal generators for $\Lambda_{\gen}$, as well as an expression for $c(\Lambda_{\gen})$ (see Corollary \ref{cond}).
\end{rmk} 

\begin{rmk} \label{1modp} If $m \not\equiv -1$ (mod $p$) and $m\neq p+1$ we have 
\begin{equation}
\label{eqn:exmZ}
g_1=p+m+1.
\end{equation}
\end{rmk} 
Indeed since $g_{-1}=p$ and $g_0=m$, then $u_0=m+p,$ for certainly $m+p\in (\Gamma+p)\cap(\Gamma+m),$ and any $\alpha\in \Gamma$ with $\alpha<p+m$ is either a multiple of $p$ or is equal to $m$. We claim $p+m+1 \not\in\Gamma.$ Otherwise there exist $a,b\geq0$ such that $(a-1)p+(b-1)m=1.$ One of these coefficients must be negative, so either $a=0$ or $b=0$. In the first case $m = p+1$ and in the second case $m = (a-1)p-1$.

In the excluded cases we have that $g_1=p+m+2.$ In either case $g_1-p$ is the {\em Zariski invariant} (see {\em e.g.} \cite{HH}). 

\begin{empl}
\label{baby}
Consider $\Gamma = \langle 10, 23 \rangle$. Since $23 = 2\cdot10 + 3$, we see by the last remark that $u_0 = m + p = 33$ and $g_1 = u_0 + 1 = 34.$ To find $u_1 = g_1 + \gamma$, write out the elements of $\Gamma$ and observe that $\gamma = 2m = 46$ is the smallest element of $\Gamma$ that satisfies the condition $g_1 + \gamma \in E_0 = \Gamma^*$, whence $u_1= 80$. Then to find $g_2 \in (\mathbb{N} + u_1) \setminus E_1$, write out the first few elements of $\Gamma + g_1$ to check that $u_1 + 1$ is {\em not} in $E_1$, so $g_2 = u_1 + 1 = 81$. Find $u_2 = g_2 + \gamma_2$ by seeing that $g_2 + m = 104 = g_1 + 7p \in E_1$ via direct inspection of $E_1$, and no element $\gamma$ smaller than $m$ exists in $\Gamma^*$ such that $g_2 + \gamma \in E_1$. One can continue in this way, but the checking becomes harder as the sets $E_i$ become more complicated.
\end{empl}

The calculation of generic generators is simplified by the following result from \cite{D2}.
\begin{lemma}
\label{lemma:Delorme}
Setting $\bar{u}_i= u_i+\mu-pm$ for $i\in (0, n)$, there exists a number $c_i\in \mathbb{Z}$ such that
\[(\mathbb{N}+\bar{u}_i)\cap E_i=(\mathbb{N}+\bar{u}_i)\cap (\Gamma+c_i).\]
\end{lemma}
\qed

\noindent
Delorme provides a recursive formula for $c_i.$ Define $c_0=0$. Then 
\begin{equation}
\label{eqn:ci}
c_{i+1}=c_i+g_{i+1}-u_{i+1}.
\end{equation}
Note that $\mu-mn<0$, so $\bar{u}_i<u_i$. Finally we note that the Delorme algorithm stops when $u_i-c_i\geq \mu$.
\begin{empl} (continued from Example \ref{baby}.)  For $\Gamma = \langle 10, 23 \rangle$ we found $u_0 = 33$, $g_1= 34$, $u_1 = 80$, and $g_2 = 81$. In particular, $u_1 = g_1+ 2m,$ so by Delorme's formula $c_1 = c_0 + (g_1 - u_1) = -2m = - 46.$ To find the smallest $\gamma \in \Gamma$ such that $u_2 = g_2 + \gamma \in E_1$, we can instead test: 
\begin{equation}
\label{test}
\text{For which $\gamma$ is \;} g_2 + \gamma \in (\Gamma+c_1)? \text{\; -- that is, \;} g_2 + \gamma-c_1\in\Gamma?
\end{equation}
Since $g_2 - c_1 + \gamma = 127 + \gamma,$ we see that $\gamma = m$, because it verifies the condition in (\ref{test}) and the only smaller values of $\gamma$ are $p$ and $2p$, which both fail condition (\ref{test}). Thus $u_2 = g_2 + m = 104$. Then to compute $g_3 = u_2 + r$ we replace the algorithm's requirement that $r$ be the minimal positive integer such that $u_2 + r \notin E_2$, with the condition $u_2 - c_2 + r \notin \Gamma$. Since $u_2 - c_2 =  u_2 + m = 127$, and $128 \notin \Gamma,$ we find that $r = 1$ and $g_3 = u_2 + 1 = 105$. One more cycle of calculation completes the algorithm: We now have $c_2 = c_1 + (g_2 - u_2) = -3m = -69$, and we seek $\gamma_3$ such that $g_3 - c_2 + \gamma_3 \in \Gamma$. Thus $\gamma_3$ is the smallest element of $\Gamma$ such that $174 + \gamma_3 \in \Gamma$. Since $8m = 184 = 174 + p$, we have $\gamma_3 = p$. So $u_3 = g_3 + \gamma_3 = 115$. To find $g_4$, we seek $r$ such that $u_3 + r - c_3 \notin \Gamma$. Since $u_3 - c_3 = 115 + 79 = 194$, we check that  $194 + r \in \Gamma$ for $r = 1, 2$ and $194+3 \notin \Gamma$. Thus $g_4 = u_3 + 3 = 118$. The algorithm now ceases: $g_4 - c_3 + \gamma_4 = 118 + 79 + \gamma_4 = 197 + \gamma_4$ is in $\Gamma$ if $\gamma_4 = p = 10$. But then $u_4 = g_4 + \gamma_4 = 128$ and $u_4 - c_4 = 217$ is greater than $\mu = 198$, so no $r > 0$ exists satisfying $u_4-c_4 + r \notin \Gamma$.  
\end{empl}

\section{Explicit calculation of the generators: preliminaries} \label{prelims}
We introduce the ingredients of our main result, which presents
explicit formulas for the generators of $\Lambda_{\gen}$.
Central to our calculation is the data provided by the Euclidean algorithm applied to $m$ and $p$. Let $s$ be the number of steps in the Euclidean algorithm for $m$ and $p$, define $p_0 = p$, and
\begin{align}
\label{eqn:ki}
m & = k_0p_0 + p_1\\
p_0 &= k_1p_1 + p_2  \nonumber \\ 
p_1 &= k_2p_2 + p_3 \nonumber\\
&\phantom{p=\,} \vdots  \nonumber\\
p_{s-2} &= k_{s-1}p_{s-1} + 1 \nonumber
\end{align}   
where $1 < p_i < p_{i-1}$ for $1 \leq i \leq s-1$. In accordance with the above, set $p_s = 1$. The number $s$ is called the \emph{level} of the semigroup. The numbers $p_i$ are the {\em divisors} and $k_i$ are the \emph{quotients} for the semigroup. Sometimes we refer to $p_{s-1}$ as the \emph{final divisor}.  It is natural to define $k_s = p_{s-1}$ and $p_{s+1} = 0$, so we may conveniently write $p_{s-1} = k_sp_s + p_{s+1}$. 

From this sequence we derive related numbers:
Let
\begin{equation}
\begin{pmatrix}
    A_0 & A_1  \\
    B_0 & B_1
  \end{pmatrix}
=  \begin{pmatrix}
    0 & 1  \\
    1 & k_0 
  \end{pmatrix}
    \label{eq:myeqn}
\end{equation}
and for $1 \leq i < s$ define 
\begin{equation}
\begin{pmatrix}
    A_{i+1} \\
    B_{i+1}
  \end{pmatrix}
=  \begin{pmatrix}
  A_{i-1} & A_i \\
    B_{i-1} & B_i
  \end{pmatrix}
\begin{pmatrix}
    1  \\
    k_i 
  \end{pmatrix}.
      \label{eqn:Ai}
\end{equation}
Then by induction we have  
\begin{align}
\label{eqn:pi}
 p_i = (-1)^{i}B_ip +  (-1)^{i-1}A_im, \;\;\; 1 \leq i \leq s.
\end{align}
In particular for $i = s$ we have
\begin{align}
\label{eqn:ps}
1 = (-1)^{s}B_sp +  (-1)^{s-1}A_sm.
\end{align}

The following simple induction will also be useful. 
\begin{lemma} 
For $1 \leq j \leq s$, 
\begin{align}
\label{eqn:p}
p & = A_jp_{j-1} + A_{j-1}p_j\\
m & = B_jp_{j-1} + B_{j-1}p_j
\label{eqn:m}
\end{align}
\end{lemma} 
\qed

The case $i = s$ in the last lemma is particularly useful:
\begin{equation}
\label{eqn:pp}
p = A_sp_{s-1} + A_{s-1} \text{\;\;and \;\;} m = B_sp_{s-1} + B_{s-1}.
\end{equation}
Recalling that the final divisor $p_{s-1}$ is also $k_s$, the last equations can be seen as an extension of (\ref {eqn:Ai}) if we take $A_{s+1} = p$ and $B_{s+1} = m$. The equations (\ref {eqn:pi}), (\ref {eqn:Ai}),  and (\ref {eqn:pp}) are central to the proof of the main theorem.

Following Delorme, we can represent generators in a standard form: Any element of $\Lambda \setminus \Gamma^*$ can be written uniquely as
\begin{equation}
\label{Delm}
g = pm - am - bp  
\end{equation}
where $ 0 < a <  p$ and $0 < b < m$. 

In addition to the numbers already defined, we set $n_s := p_{s-1}$, and for $1 \leq l \leq s,$ denote $N_l=\sum_{j=l}^s n_{j},$ and $N_{s+1}=n_{s+1}=0$, where we recursively define
\begin{equation}
\label{eqn:ns}
n_l =
    \begin{cases}
    0, & \text{if } 2\mid N_{l+1} \text{ and } n_{l+1} \neq 0, \\
        k_l, & \text{if }  2 \nmid N_{l+1} \text{ or } n_{l+1}=0.\\
    \end{cases}
\end{equation}
Note that $n_s$ is never zero, and it is impossible for two consecutive values of $n_l$ to be zero. 

\section{Explicit calculation of the generators: Main Theorem} \label{MT}

As always $\Gamma$ is the semigroup generated by coprime integers
$p$ and $m$. The case where $p=2$ is trivial, and henceforth we assume
$p>2$. 

\begin{theorem} \label{Main}  For $1\leq i\leq N_1-1$ let
\begin{align*}
u_i &= g_i + \gamma_i, \;\text{and} \\
g_{i+1} &=u_i + p_{j}
\end{align*}
where $j \in [1, s]$ is determined by $N_{j+1}\leq i\leq N_{j}-1$ and
\begin{equation}
\gamma_i=
\begin{cases}
(B_{j} - 1)p & \text{ if } \quad 2\nmid j \text{\;\;and\;\;} 2\nmid i, \\
\;\;\;\;\;\;\;p & \text{ if } \quad 2\nmid j \text{\;\;and\;\;} 2\mid i, \\
(A_{j} - 1)m & \text{ if } \quad 2\mid j \text{\;\;and\;\;} 2\nmid i, \\
\;\;\;\;\;\;\;m & \text{ if } \quad 2\mid j \text{\;\;and\;\;} 2\mid i.
\end{cases}
\end{equation}

Then the generic $\Gamma$-semimodule $\Lambda_{\gen}$ is generated by
$g_{-1}=p$, $g_0=m$, $g_1=p+m+1$, and by
\begin{equation} \label{recur}
g_{i+1}=g_i + \gamma_i + p_{j}
\end{equation}
for $1\leq i\leq N_1-1$. 
 \end{theorem}
\bigskip

Our notation agrees with that used in Delorme's algorithm, and the generators produced by recursion (\ref{recur}) are minimal except in circumstances we now explain. Some generators produced by Theorem \ref{Main} are not minimal in two situations: (a) The first inequality constraining $i$ in the statement does not build in the stopping condition: ``Stop when $u_i - c_i \geq \mu$". (b) It is possible for particular $\gamma_i$ defined in the theorem to be zero. This happens when $k_0 = 1$ or $k_1 = 1$. In Section \ref{specialcases} we identify non-minimal generators produced by the recursion in this situation. Outside of this case, the correct stopping point for $i$ and thus the precise identification of the set $G$ of minimal generators is given in Corollary \ref{cond}. The correct identification of $G$ and $|G|$ when $k_0 = 1$ or $k_1 = 1$ is given in Section \ref{specialcases}.

A general formula for the conductor of $\Lambda_{\gen}$ is given in Corollary \ref{condc}.

In Section~\ref{examples} we present worked-out examples, and in Section~\ref{Tjurinaformula} relate our calculation to that of the (minimal) Tjurina number.

\begin{proof}[Proof of Theorem \ref{Main}]
First we recapitulate the algorithm of Delorme, relying heavily on Lemma~\ref{lemma:Delorme}. Each generator is calculated in two steps: (1) From the last generator $g_i$, a ``collision" $u_i = \gamma_i + g_i $ is computed. It is the smallest value of the form $\gamma + g_i$, $\gamma \in \Gamma$, which belongs to $E_{i-1}$, the set generated by the previous generators $\{g_j\}_{j \leq i-1}$ under the action of $\Gamma$. (2) The next generator, $g_{i+1}$, is found by taking the ``minimal jump" from $u_i$; that is, $g_{i+1} = u_i + r_i$ where $r_i$ is the smallest positive integer $r$ such that $u_i + r$ is not in $E_i$. The calculations of $u_i$ and $g_{i+1}$ are simplified by Lemma \ref{lemma:Delorme}, which allows us to use $\Gamma + c_{i'}$ in the role of $E_{i'}$ ($i' = i-1$ or $i$), where $c_{i'} = c_{i'-1} +g_{i'} - u_{i'}$ and $c_0 = 0$. This reduces calculating $r_i$ to finding the least $r > 0$ satisfying $u_i- c_i + r  \not\in \Gamma$, and similarly for $\gamma_i$ finding the least $\gamma \in \Gamma$ such that $ \gamma + g_i - c_{i-1} \in \Gamma$.

We will need two levels of induction:  on the index of the generator $i$, and on the level $j$ of the Euclidean Algorithm (the latter begins at $j=s$ and decreases). The proof will focus on the numbers $g_i, u_i$ and $c_i$. It follows from the recursive property of $c_i$ (see Lemma \ref{lemma:Delorme}) and the definition of $\gamma_i$ that $c_i = -\sum_{a=1}^i \gamma_a$. We first prove the base case ($j=s$) for $1\leq i\leq k_s-1$.  
\\
\\
{\bf Base Case $j=s$:} We assume $2\nmid s$, since the calculations for $s =$ (even) are exactly parallel (we will however provide the equivalent intermediate expressions). By definition $c_0=0,$ and so $g_1-c_0=g_1.$ Thus we seek the smallest $\gamma_1\in\Gamma$ such that $g_1+\gamma_1\in E_0=\Gamma^*.$ This is simplified by expressing $g_1 -c_0$ in standard form (\ref{Delm}). First note that by (\ref{eqn:ps}) and $2\nmid s$
\begin{equation}
\label{ps=1}
1=A_sm-B_sp. 
\end{equation}
Because $g_1$ is given by $m + p +1$ (see (\ref{eqn:exmZ})), and $p=k_sA_s+A_{s-1}$ by (\ref{eqn:pp}), we obtain
\begin{equation}
\label{step1}
g_1 - c_0 = pm-((k_s-1)A_s+A_{s-1}-1)m-(B_s-1)p.
\end{equation}
Now compare the terms of (\ref{step1}). The inequality $k_s = p_{s-1} > 1$  ensures 
$$((k_s-1)A_s+A_{s-1}-1)m\geq A_s m>B_s p>(B_s-1)p,$$
and so the minimal element in $\Gamma$ to add to $g_1$ must be $\gamma_1=(B_s-1)p,$ as any smaller element would result in negative coefficients for both $p$ and $m.$ The algorithm sets $u_1 = g_1 + \gamma_1$, so we have $u_1=pm-((k_s-1)A_s+A_{s-1}-1)m$ and $c_1=-(B_s-1)p$. Therefore
\begin{equation}
u_1-c_1=pm-((k_s-1)A_s+A_{s-1}-1)m+(B_s-1)p\in \Gamma.
\end{equation}
We continue now to show that $r_1 = 1$. Note that $1 = p_s$, as expected in the base case. By (\ref{ps=1}) we have
\[ 
u_1-c_1+1=pm-((k_s-2)A_s+A_{s-1}-1)m-p.
\]
If the coefficient of the middle term is zero, then $u_1 - c_1 > \mu$ and the algorithm stops with $g_1$. Otherwise $u_1 - c_1 + 1$ is not in $\Gamma$, and so $r_1 =1$. Thus $g_2 = u_1 + 1$ and we have \[g_2 - c_1 = u_1 -c_1 + 1 = pm - ((k_s-2)A_s + A_{s-1} - 1)m - p.\]
This establishes the first step of induction for the base case.
\begin{rmk} \label{struc} The structure of the calculation that emerges has the following form:
\begin{enumerate}
\item Begin the step by expressing $g_i -c_{i-1} $ in standard form (\ref{Delm}): $$g_i -c_{i-1} = pm - a_im - b_im.$$
\item Calculate $\gamma_i$ as the smaller of the two terms $a_im$ and $b_ip$.
\item Because $u_i = g_i + \gamma_i$ and $c_i = c_{i-1} - \gamma_i$, it is simplest to think of $u_i - c_i$ as obtained from $g_i-c_{i-1}$ by reversing the sign of the smaller term: 
\[
u_i - c_i = \begin{cases}
pm + a_im - b_ip,  & \gamma_i = a_im < b_ip, \\
pm - a_im + b_ip, & \gamma_i = b_ip < a_im.
\end{cases}
\]
\item Find the smallest $r_i > 0$ such that $u_i - c_i + r_i \notin \Gamma$. We show below that $r_i = p_j$ for level $j$, and compute the interval of values $i$ belonging to this level. Then $g_{i+1} = u_i + r_i$ and the process repeats. 
\end{enumerate}
\end{rmk}
\noindent
{\bf Induction step for base case.} Suppose $1 \leq i < k_s$, and inductively assume
\begin{equation} \label{base}
g_{i}-c_{i-1}= pm-((k_s-i)A_s+A_{s-1}-1)m-
\begin{cases}
(B_s - 1)p & \text{ if } \quad 2\nmid i, \\
p & \text{ if } \quad 2\mid i. \\
\end{cases}
\end{equation}
If $2\mid i$ in (\ref{base}) then clearly $\gamma_{i}=p$. 
Thus 
$$u_i-c_i=pm-((k_s-i)A_s+A_{s-1}-1)m+p,$$
and, provided $i < k_s$, we check that $r_i = 1$: 
\begin{equation}
\label{level1ri}
u_i-c_i+1=pm-((k_s-i-1)A_s+A_{s-1}-1)m-(B_s-1)p\notin\Gamma
\end{equation}
since both coefficients are negative (here we assume $A_{s-1}>1$; see Section \ref{specialcases}). Now suppose $2\nmid i$. Then since $k_s - i >0$ we have  $((k_s-i)A_s+A_{s-1}-1)m>(B_s-1)p.$ Hence $\gamma_i=(B_s-1)p.$
A computation as above shows that $r_i=1$ again. Thus we have $r_i = 1$,
\begin{align*}
\gamma_i=
\begin{cases}
(B_s - 1)p & \text{ if } \quad 2\nmid i, \\
\;\;\;\;\;\;\;p & \text{ if } \quad 2\mid i
\end{cases}
\end{align*}
and $g_{i+1}-c_i$ has the same form as equation (\ref{base}). Thus the induction step is proved. In particular, for $i = k_s-1$  we find if $2 \nmid s$ then
\begin{equation}
\label{base1}
g_{N_s}-c_{N_s-1}= pm-(A_{s-1}-1)m-
\begin{cases}
(B_s - 1)p & \text{ if } \quad 2\nmid N_s, \\
p & \text{ if } \quad 2\mid N_s. \\
\end{cases}
\end{equation}
The case  $2 \mid s$ is exactly analogous with roles reversed: we now have $B_sp - A_sm =1$, which provides $B_sp > A_sm$, and so $\gamma_1 = (A_s -1 )m$, etc. Thus if $2 \mid s$ then
\begin{equation}
\label{base2}
g_{N_s}-c_{N_s-1}= pm-(B_{s-1}-1)p-
\begin{cases}
(A_s - 1)m & \text{ if } \quad 2\nmid N_s, \\
m & \text{ if } \quad 2\mid N_s. \\
\end{cases}
\end{equation}
This establishes the case where $j=s$ and $1 \leq i\leq N_s-1.$ Note that in this case we need not prove $r_i = 1$ is the {\em minimal} jump, since $r_i$ must be positive. Lastly, if $A_{s-1} = 1$, the generator $g_{N_s}$ belongs to $\Gamma + c_{N_s-1}$, and so is not minimal; see Section \ref{specialcases}.  

{\bf Induction on level $j < s$.} We assume the theorem formulas for $j+1$ and prove them for $j$. We illustrate the pattern of moving from an odd level to an even level, hence assume $2\mid j$. The other case will be clear with obvious reversal of roles. 
We will also first assume $n_{j+1}\neq 0$.

Assuming $2\mid j,$ and $n_{j+1}\neq 0,$ we have by induction (compare (\ref{base1}))
\begin{equation}
\label{ind1}
g_{N_{j+1}}-c_{N_{j+1}-1}= pm-(A_{j}-1)m -
\begin{cases}
(B_{j+1}-1)p & \text{ if } 2\nmid N_{j+1} \\
p & \text{ if } 2\mid N_{j+1} 
\end{cases}
\end{equation}

Here for reference  is the equivalent induction statement  assuming $2\nmid j$ (compare (\ref{base2})):
\begin{equation}
\label{ind2}
g_{N_{j+1}}-c_{N_{j+1}-1}= pm-(B_{j}-1)p -
\begin{cases}
(A_{j+1}-1)m & \text{ if } 2\nmid N_{j+1} \\
m & \text{ if } 2\mid N_{j+1} 
\end{cases}
\end{equation}

We first consider the case where $2\nmid N_{j+1}.$ Since $2\mid j$, by (\ref{eqn:pi}) we have 
\begin{equation} 
\label{eqn:pj}
 p_j = B_jp - A_jm.
\end{equation}
This implies $(A_{j}-1)m<(B_{j+1}-1)p$, since we always have $p<m$ and $B_j < B_{j+1}$. Thus (see Remark \ref{struc}) we have $\gamma_{N_{j+1}}=(A_{j}-1)m.$

It follows that 
\begin{equation} \label {delta}
u_{N_{j+1}}-c_{N_{j+1}}=pm+(A_{j}-1)m-(B_{j+1}-1)p.
\end{equation}
We claim $r_{N_{j+1}}$ is $p_j.$ Indeed by (\ref{eqn:pj})
\begin{equation}
\label{uplus}
u_{N_{j+1}}-c_{N_{j+1}}+p_{j}=pm-m-(B_{j+1}-B_{j}-1)p \notin\Gamma
\end{equation}
since $B_{j+1} > B_{j}+1$. In order to show that $p_{j}$ is minimal, we first relate any positive number $r$ to the divisors $p_i$ in the Euclidean algorithm (\ref{eqn:ki}).

\begin{defn}  \label{rrecur} Suppose $0 < r < p_j$. Set $r^1 = r$, and for $a > 0 $ define $\{r^{1 + a}, \alpha_{j+a}\}$ by the conditions $r^{a} = \alpha_{j+a}p_{j+a} + r^{1 + a}$ and $0 \leq r^{1+ a} < p_{j+a}$. Let $h = \text{min}\{j + a | r^{1+a} = 0\}$.
\end{defn}
Returning to the proof, suppose $r < p_j$. Then with $\{\alpha_i\}$ as in Definition \ref{rrecur} we have $r= \sum_{i=j+1}^h \alpha_ip_i$. By (\ref{eqn:pi}) we may write $p_i=(-1)^{i+1}(A_im-B_ip)$, and 
\setlength{\jot}{8pt}
\begin{align} \nonumber
r  	& = \sum_{i=j+1}^h(-1)^{i+1}\alpha_i(A_im-B_ip), \\
	& = \left(\sum_{i=j+1}^h(-1)^{i+1}\alpha_i A_i \right)m +  \left(\sum_{i=j+1}^h(-1)^i\alpha_i B_i\right)p. \label{rsums}
\end{align}
We will apply the following lemma to the analysis of $r$.

\begin{lemma} \label{sumpos} 
Let  $x = \sum_{i=u}^h (-1)^{i+1} \alpha_i A_i$ and $y = \sum_{i=u}^h (-1)^i \alpha_i B_i$, where $u \leq h \leq s$ and $0 \leq \alpha_i \leq k_i$. Assume $\alpha_h > 0$.
\begin{enumerate}
\item $x > 0 \Leftrightarrow 2 \nmid h$, and $y > 0 \Leftrightarrow 2 \mid h$. 
\item If in addition $\alpha_i = k_i \Rightarrow \alpha_{i+1} = 0$ then 
\begin{enumerate}
\item $x > 0 \Rightarrow x \geq A_u$, and $y > 0 \Rightarrow y \geq B_u$.
\item $x < 0 \Rightarrow x > A_s -p$, and $y < 0 \Rightarrow y > B_s - m$.
\end{enumerate}
\end{enumerate}
\end{lemma}
We defer the proof of Lemma \ref{sumpos} to the end of this section, and proceed with the proof of the main theorem. Write $r = xm + yp$, where $x$ and $y$ are the sums in (\ref{rsums}). By the equations (\ref{eqn:ki}) for $p_i$, and since $r < p_j$, the recursion in Definition \ref{rrecur} ensures $\alpha_{i} \leq k_i$ and $\alpha_{i} < k_i$ unless $\alpha_{i+1} = 0$. So $x$ and $y$ satisfy Part 2 of Lemma \ref{sumpos}, with $u = j+1$. 

Now consider $u_{N_{j+1}}-c_{N_{j+1}}+r$ with $r < p_j$. From (\ref{delta})
\begin{equation}
u_{N_{j+1}}-c_{N_{j+1}} + r = pm+(A_{j} + x -1)m - (B_{j+1} - y -1)p. \\
\end{equation}
Either $x > 0$ and $y < 0$, or the reverse (since $r < p$). Suppose first $x > 0$ and $y < 0$. Then 
\[ u_{N_{j+1}}-c_{N_{j+1}} + r = (A_{j} + x -1)m + (m - B_{j+1} + y +1)p \]
Then by Lemma \ref{sumpos} we have $y > B_s - m$, and both coefficients are positive. Now suppose $x < 0$ and $y > 0$. Then 
\[ u_{N_{j+1}}-c_{N_{j+1}} + r = (p + A_{j} + x -1)m + (y - B_{j+1} +1)p \]
where Lemma \ref{sumpos} gives $y \geq B_{j+1}$ and $x > A_s - p$. Thus once again both coefficients are positive, so in either case adding $r$ to $u_{N_{j+1}}-c_{N_{j+1}}$ results in an element of $\Gamma$. Therefore with the starting assumption $2\nmid N_{j+1}$, we see $r = p_j$ is the minimal jump needed to escape $\Gamma$. Therefore $r_{N_{j+1}} = p_j$ and (compare (\ref{uplus}))
\begin{equation}
g_{N_{j+1}+1}-c_{N_{j+1}}= pm-m-(B_{j+1}-B_{j}-1)p.
\end{equation}
The pattern is now set for $k_j$ steps. For $0 \leq i \leq k_j - 1$ we have $r_{N_{j+1}+i} = p_j$,
\begin{align*}
\gamma_{N_{j+1} + i}=
\begin{cases}
(A_j - 1)m & \text{ if } \quad 2\nmid (N_{j+1} + i), \\
\;\;\;\;\;\;\;m & \text{ if } \quad 2\mid (N_{j+1} + i)
\end{cases}
\end{align*}
and $g_{N_{j+1}+i+1} = u_{N_{j+1}+i} + \gamma_{N_{j+1}+i}$, so
\begin{equation}
g_{N_{j+1} +1+ i}-c_{N_{j+1} + i}= pm-(B_{j+1} - iB_j-1)p-
\begin{cases}
(A_j - 1)m & \text{ if } \quad 2\nmid (N_{j+1} + i) \\
m & \text{ if } \quad 2\mid (N_{j+1} + i). \\
\end{cases}
\end{equation}

Now suppose instead we have that $2\mid N_{j+1}.$ Then the induction statement (\ref{ind1}) is
\begin{equation}
\label{Neven}
g_{N_{j+1}}-c_{N_{j+1}-1}=pm-(A_{j}-1)m-p,
\end{equation}
and we can assume $A_j > 1$ (otherwise the algorithm ends). Hence in this case $\gamma_{N_{j+1}}=p$,
and we have 
$$u_{N_{j+1}}-c_{N_{j+1}}=pm-(A_{j}-1)m+p.$$
First note that $r_{N_{j+1}} \neq p_j$, since $u_{N_{j+1}}-c_{N_{j+1}}+p_{j}=(p-2A_{j}+1)m+(B_{j}+1)p\in\Gamma$.
However
\[u_{N_{j+1}}-c_{N_{j+1}}+p_{j-1}=pm-(A_{j}-A_{j-1}-1)m-(B_{j-1}-1)p\notin\Gamma\]
since otherwise either $k_1 = 1$ or $A_j = 2$ (in the latter case $u_{N_{j+1}} - c_{N_{j+1}}$ exceeds $\mu$). One shows that $p_{j-1}$ is the minimal value by analyzing $r < p_{j-1}$ via Lemma \ref{sumpos}, exactly as in the last argument. It follows that in this case
$r_{N_{j+1}}=p_{j-1}$ and $g_{N_{j+1}+1} = u_{N_{j+1}} + p_{j-1}$, so
\begin{equation}
\label{gNj}
g_{N_{j+1}+1}-c_{N_{j+1}}=pm-(A_{j}-A_{j-1}-1)m-(B_{j-1}-1)p.
\end{equation}
Because $A_j = k_{j-1}A_{j-1} + A_{j-2}$ we see again by (\ref{eqn:pi}) that $(B_{j-1} - 1)p$ is the smallest element of $\Gamma$ which added to (\ref{gNj}) results in an element of $\Gamma$. The previous pattern now repeats: $\gamma_{N_{j+1}+i}$ alternates between $(B_{j-1}-1)p$ and $p$, while $r_i = p_{j-1} = A_{j-1}m - B_{j-1}p$ causes the coefficient of $m$ to get closer to zero. Thus for $k_{j-1}$ steps we have
\begin{equation}
g_{N_{j+1} + i + 1}-c_{N_{j+1} + i}= pm-(A_j - iA_{j-1}-1)m-
\begin{cases}
(B_{j-1} - 1)p & \text{ if } \quad 2\nmid (N_{j+1} + i), \\
p & \text{ if } \quad 2\mid (N_{j+1} + i). \\
\end{cases}
\end{equation}
The assumption $(2 \mid N_{j+1})$ leads to calculations of $\gamma_i$ and $r_i$ corresponding to level $(j-1)$ instead of $j$ for an interval of length $k_{j-1}$, which justifies the definition (\ref{eqn:ns}) of $n_j = 0$ and $n_{j-1} = k_{j-1}$ in this case. 

Now supposing that $n_{j+1} = 0$, then by (\ref{eqn:ns}) we have  $N_{j+1} = N_{j+2}$, and $n_{j+2} \not = 0$. The inductive assumption is given by (\ref{ind2}), with $j$ replaced by $(j+1)$:
\begin{align*}
g_{N_{j+1}}-c_{N_{j+1}-1} &= g_{N_{j+2}}-c_{N_{j+2}-1} \\
					&= pm-(B_{j+1}-1)p - m.
\end{align*}
With the same analysis used in (\ref{Neven}), we have $\gamma_{N_{j+1}} = m$, and $r_i = p_j$. Thus level $(j+1)$ is empty, and level $j$ proceeds for $k_j$ steps, justifying the definition $n_j = k_j$ in this case. Lastly, the analysis for $2 \nmid j$ is clearly the strict analog of the cases just presented.
\end{proof}

We now prove Lemma \ref{sumpos}.
\begin{proof} We prove all the statements for $x$ only, since the arguments for $y$ are strictly analogous.

{\bf Proof of 1 ($\Leftarrow$).} Suppose $2\nmid h$ and proceed by induction. First suppose $h = u$. Then 
\[ x = \alpha_u A_u > 0. \] Now suppose the implication is true for $h' < h$. Then
\[x = \alpha_h A_h - \alpha_{h-1}A_{h-1} + \sum_{i=u}^{h-2} (-1)^{i+1} \alpha_i A_i  \]
(The sum is empty if $h = u+1$.) By assumption $\alpha_h A_h - \alpha_{h-1} A_{h-1} \geq A_h - \alpha_{h-1} A_{h-1} = (k_{h-1} - \alpha_{h-1}) A_{h-1} + A_{h-2}$. Thus
\begin{equation}
\label{compr}
 x \geq (k_{h-1} - \alpha_{h-1}) A_{h-1} + \left ( A_{h-2} + \sum_{i=1}^{h-2} (-1)^{i+1} \alpha_i A_i \right )
\end{equation}
The first expression is non-negative since $\alpha_{h-1} \leq k_{h-1}$. Adding $A_{h-2}$ to the last sum guarantees the coefficient of the $(h-2)$ term, namely $\alpha_{h-2} + 1$, is positive, so the combined sum is also positive by induction. (If $\alpha_{h-2} + 1$ exceeds $k_{h-2}$, the violation of the hypothesis is in the positive direction.)

{\bf Proof of 1 ($\Rightarrow$).} Suppose $2 \mid h$. Apply the previous argument to $-x$, as the only use made of the parity of $h$ is that the highest-index term has positive coefficient.

{\bf Proof of 2(a).} If $h = u$ and $x > 0$ then $x = \alpha_uA_u \geq A_u$. Suppose $h > u$. Then (\ref{compr}) is true and the last sum is positive. Because $\alpha_h > 0$ by assumption, the additional condition ($\alpha_i = k_i \Rightarrow \alpha_{i+1} = 0$) implies that $(k_{h-1} - \alpha_{h-1})$ is positive, so $x > A_{h-1} \geq A_u$. 

{\bf Proof of 2(b).} Assume $x < 0$. First suppose $h = u$. Then $2 \mid u$ and $x = -\alpha_uA_u$. If $u < s$ then 
\[ x \geq -k_uA_u = -(A_{u+1} - A_{u-1}) >  -A_s  >  -(k_s-1)A_s - A_{s-1} = A_s - p \] 
(see (\ref{eqn:pp}) and recall $k_s = p_{s-1} > 1$). If instead $u = s$, then by assumption $\alpha_ u = \alpha_s \leq k_s -1$. Then 
\[x \geq -(k_s-1)A_s > -(k_s-1)A_s - A_{s-1} = A_s - p.\]
Now suppose $h > u$ (and $2 \mid h$). Then 
\[x = -\alpha_hA_h - A_{h-1} + \left( A_{h-1} +\sum_{i=u}^{h-1} (-1)^{i+1} \alpha_i A_i \right) \]
The parenthetical term is positive by Part 1. Therefore 
\[x > -\alpha_hA_h - A_{h-1} \geq -k_hA_h - A_{h-1} = -A_{h+1}. \]
If $h < s$ then $-A_{h+1} \geq -A_s > -(k_s-1)A_s - A_{s-1} = A_s - p$. If $h = s$ then $\alpha_h = \alpha_s \leq k_s-1$ and
\[ x > -(k_s-1)A_s - A_{s-1} = A_s -p. \]
\end{proof}
Having established Lemma \ref{sumpos}, the proof of the main theorem is now complete. We next refine the main result with two corollaries.

\begin{corollary} \label{condc}
Let $n$ be the index of the last minimal generator $g_n$ produced by the recursion in Theorem \ref{Main}. The conductor of $\Lambda_{\gen}$ is given by
\begin{equation*}
c(\Lambda_{\gen}) = \mu + c_n.
\end{equation*}
\end{corollary}
\begin{rmk} The index $n$ corresponds to the last {\em minimal} generator produced by Theorem \ref{Main}, which may occur before the last value of $i$ is reached. It's value is determined in Corollary \ref{cond} and Section \ref{specialcases}.
\end{rmk}
\begin{proof}
By Lemma \ref{lemma:Delorme} we have 
\begin{equation} \label{lemcon}
(\mathbb{N}+\bar{u}_n)\cap E_n=(\mathbb{N}+\bar{u}_n)\cap(\Gamma+c_n).
\end{equation}
Recall that $E_n = \Lambda_{\gen}$, and notice the obvious fact that the conductor of $(\Gamma + c_n)$ is just $\mu + c_n$. Part (ii) now follows from

{\bf Claim:} $\bar{u}_n < \mu + c_n$.

For if the Claim is true then (\ref{lemcon}) implies that $\Lambda_{\gen}$ and $\Gamma + c_n$ have the same conductor. But $g_n - c_{n-1}$ is of the form $pm - \alpha m - \beta p$ where $\alpha, \beta >0$, and $\gamma_n = \text{min}\{\alpha m, \beta p \}$. Clearly $\alpha < p$ and $\beta < m$, so $\alpha m \neq \beta p$. It follows that $u_n - c_n = g_n - c_{n-1} + 2\gamma_n < pm$. From Lemma \ref{lemma:Delorme} we have $\bar{u}_n=u_n + \mu - pm,$ therefore $\bar{u}_n - c_n = u_n - c_n + \mu - pm < pm + \mu - pm = \mu$.
\end{proof}

\begin{rmk} Note that $c_n = -\sum_{j=1}^n \gamma_j$, so the conductor can be calculated from the recursion of Theorem \ref{Main}. 
\end{rmk}
\begin{rmk} Except in the extreme case where either $k_0$ and/or $k_1$ equals 1 (discussed in the next section), we have  $g_n - c_n = pm - m$, so typically  $g_n - (p-1)$ equals the conductor of $\Lambda_{\gen}$.
\end{rmk}

\begin{corollary} \label{cond} Suppose $k_0, k_1\neq 1$, and let $g_n$ be the last minimal generator produced by Theorem \ref{Main}. Let $G$ be the set of minimal generators for $\Lambda_{\gen}$. Then  $G = \{g_i| -1 \leq i \leq n\}$ with $g_i$ as in  (\ref{recur}), $|G|=n+2$ its cardinality, and $n$ is either $N_1$ or $N_1-2$. More precisely, the cardinality of $G$ is given by
\begin{equation*}
|G|=
\begin{cases}
N_1 + 2&\ \text{if } n_1 = 0, \\ 
N_1 &\ \text{otherwise}.
\end{cases}
\end{equation*}
\end{corollary}

\begin{proof} Theorem \ref{Main} establishes that the recursion (\ref{recur}) aligns with Delorme's algorithm. It remains to show when the recursion should end.  There are three cases for how the last level, $j=1$, can begin. If $n_2 \not = 0$ then taking the inductive statement (\ref{ind2}) with $j=1$ we have
\begin{equation}
\label{eqn:cases}
g_{N_2}-c_{N_2-1} = pm - (k_0-1)p -
\begin{cases}
m &\ \text{if } 2\mid N_2, \\
(k_1-1)m &\ \text{if } 2\nmid N_2.
\end{cases}
\end{equation}
If however $n_2 = 0$ then we saw that level 2 is empty. In this case the start of level 1 is given by (\ref{ind1}) with $j=2$  and $N_3 = N_2$:
\begin{equation} \label{last2}
g_{N_2}-c_{N_2-1}=pm-(k_1-1)m-p, \text{\;\;\;where\;\;} 2\mid N_2, \text{\;\;and\;\;} n_2 = 0.
\end{equation}

In the first case of (\ref{eqn:cases}), $n_1 = 0$ by (\ref{eqn:ns}) and $\gamma_{N_2} = (k_0-1)p<m$. Therefore $u_{N_2}-c_{N_2}=pm-m + (k_0 - 1)p$ exceeds the conductor $\mu$, since we assume $k_0 \not = 1$. So $g_{N_2}$ is the final generator. Since $n_1=0$ we have $n = N_2 = N_1$ and $|G| = N_1 + 2$. 

In the other two cases, we have $n_1 = k_1$. Then $r_{N_2+i} = p_1$, and $\gamma_{N_2+i}$ alternates between $p$ and $(k_0-1)p$. First suppose $k_1 > 2$. After $k_1 -2$ steps we have
\begin{equation}
g_{N_1-2}-c_{N_1-3}= pm-m-
\begin{cases}
p & \text{ if } \quad 2\mid N_1 \\
(k_0 - 1)p & \text{ if } \quad 2\nmid N_1 \\
\end{cases}
\end{equation}
and $g_{N_1-2}-c_{N_1-2}=pm-m$. Thus $g_{N_1-2}$ is the final generator. We have $n = N_1 - 2$ and $|G| = n+ 2 = N_1$.

If on the other hand $k_1=2,$ we have in both cases $g_{N_2}-c_{N_2}=pm-m$. So $g_{N_2}$ is the final generator, and we have $n = N_2$ and $|G| = N_2 + 2$. But now $N_1 = N_2 + k_1 = N_2 + 2$, so again $|G| = N_1$.

\end{proof}

\begin{rmk} The recursion of Theorem \ref{Main} may stop before $u_i - c_i \geq \mu$ is satisfied. Indeed if $n_1 = 0$ then the final $u_i$  occurs at $i = N_1-1$ even though $u_{N_1-1} - c_{N_1-1} < \mu$. In this case the theorem allows the next generator to be defined, namely $g_{N_1}$. Then Corollary \ref{cond} shows that the next output in Delorme's algorithm, {\em i.e.}  $u_{N_1} - c_{N_1}$, does exceed $\mu$, and so the algorithm also stops, and $g_{N_1}$ is the last minimal generator.
\end{rmk} 

\section{Non-minimal generators} \label{specialcases}

In Delorme's algorithm, $\gamma_i$ is the least element of $\Gamma$ such that $u_i = g_i + \gamma_i \in E_{i-1}$. Thus non-minimal generators $g_i$, $i \leq n$, arise in the recursion of Theorem \ref{Main} iff $\gamma_i=0.$ This occurs when $A_j=1$ or $B_j=1.$ We always have $A_1=1$ and $B_0=1$, but the index for $A_j$ in any $\gamma_i$ is even, and the index for any $B_j$ is odd. It is possible however to have $A_2=1$ or $B_1=1.$ This is equivalent to the cases (a) $k_1=1$ or (b) $k_0=1$ respectively. In these cases non-minimal generators may arise in the recursion before we have reached the conductor $c(\Lambda_{\gen})$. Thus 
\[G =\{g_i| -1 \leq i \leq n\ \text{and\;} \gamma_i \neq 0 \},\] 
and the cardinality of $G$ is decreased by the number of times $\gamma_i = 0$.  
 
We summarize the effect on $|G|$ of the various configurations of extreme $k_0, k_1$. (The list of non-minimal generators in the last column can be empty.)


  \begin{center}
  \begin{table}[h]
{\tabulinesep=1.2mm
   \begin{tabu} {c|l|c|cl}
\multicolumn{2}{c|}{Constraints} & $|G|$ & Non-minimal generators  \\  \hline
      $k_0 = 1, k_1 > 1$ &$2 \mid N_2$ & $N_1 -  \floor{\frac{n_1-1}{2}}$ & $\{g_{N_2+2j-1} | 1 \leq j \leq  \floor{\frac{n_1-2}{2}}$  \} \\      \cline{2-4}     
      & $2 \nmid N_2$ & $N_1 -  \floor{\frac{n_1}{2}}$  & $\{g_{N_2+2j} | 0 \leq j \leq  \floor{\frac{n_1-3}{2}}$  \}   \\\hline
      $k_0 = 1, k_1 = 1$ &$n_3 = 0$ & $N_1 -  \floor{\frac{n_2-1}{2}}$ & $\{g_{N_3+2j-1} | 1 \leq j \leq  \floor{\frac{n_2-1}{2}}$  \} \\             \cline{2-4} 
      & $n_3 \neq 0$ & $N_1 -  \floor{\frac{n_2}{2}}$  & $\{g_{N_3+2j} | 0 \leq j \leq  \floor{\frac{n_2-2}{2}}$  \}   \\\hline
      $k_0 > 1, k_1 = 1$ &$n_1 = 0$ & $N_1 -  \floor{\frac{n_2-2}{2}}$ & If $2 \mid N_3:$ $\{g_{N_3+2j-1} | 1 \leq j \leq  \floor{\frac{n_2-1}{2}}$  \}, \\       \cline{2-3}       
      & $n_1 \neq 0$ & $N_1 -  \floor{\frac{n_2}{2}}$  & otherwise: $\{g_{N_3+2j} | 0 \leq j \leq  \floor{\frac{n_2-2}{2}}$  \}   \\\hline
        \end{tabu}}
\caption{Non-minimal generators}\label{nonmin}
  \end{table}
  \end{center}

\begin{rmk}
In the extreme cases treated in this section,  the value of the index $n$ of the last minimal generator given by the recursion of Theorem \ref{Main} can be deduced from Table \ref{nonmin}. It is always one more than the index of the last non-minimal generator if that set is non-empty, or $|G|-2$ otherwise. Equivalently, $n = |G|-2 +^\#\{\text{non-minimal generators}\}$.

\end{rmk}

   
\section{Examples} \label{examples}
\begin{empl} Recall Example \ref{baby} with semigroup  $\Gamma = \langle 10, 23 \rangle$. We compute 
\begin{align*}
23 & = 2 \cdot 10 + 3 \nonumber \\
10 & = 3 \cdot 3 + 1 \nonumber 
\end{align*}   
Thus the level is $s = 2$, and we easily compute the following table.
\begin{table}[H]
  \begin{center}
    \begin{tabular}{c | c | c | c  | c | c | c  }
          $i$ & $p_i$ & $k_i$ & $n_i$ &  $N_i$ & $A_i$ &  $B_i$  \\
      \hline
      0  & 10 & 2 & -  &  -  & 0  & 1  \\
      1  & 3 & 3 & 3  &  6 &  1 &  2 \\
      2  &  1 & 3 & 3  &  3 &  2 &  7 
    \end{tabular}
  \end{center}
\end{table}

Here $n_1 \neq 0$, so by Corollary \ref{cond} there are $N_1 = 6$ generators, including $p = 10$ and $m = 23$ and $g_1 = 1 + p + m = 34$. Following Theorem \ref{Main}, we compute each $\gamma_i$ with $1 \leq i \leq N_1 - 1 = 5$ and corresponding level and jump, and resulting $g_i$ and $u_i$ (while Corollary \ref{cond} indicates we should stop at $i = 4$, since $n = N_1 - 2$ in this case).
\begin{table}[H]
  \begin{center}
    \label{tab:table1}
    \begin{tabular}{c | c | r c | c | c | c  }
          $i$ & $j$ & $\gamma_i$ & &   $r_i=p_j$ & $g_i$ &  $u_i$  \\
      \hline
      1  & 2 & $(A_2-1)m \;=$ 	& 46 & 1 & 34 & 80 \\
      2  & 2 & $m \;=$ 			& 23 & 1  & 81 & 104  \\
      3  & 1 & $(B_1-1)p \;=$ 	& 10 & 3 & 105 & 115  \\
      4  & 1 & $p \;=$ 			& 10 & 3 & 118 & 128 \\
      5  & 1 & $(B_1-1)p \;=$ 	& 10 & 3 & (131) & (141)  \\
    \end{tabular}
  \end{center}
\end{table}
\end{empl} 
Note that $c_4 = -\sum_{a=1}^4 \gamma_a = -89$, so $u_4 - c_4 = 217$, which is greater than the conductor of $\Gamma$. So the algorithm stops at $i=4$ and the last displayed generator is redundant, as Corollary \ref{cond} implies it should be. Lastly we find the conductor of $\Lambda_{\gen}$:
\[ c(\Lambda_{\gen}) = \mu + c_n = (23-1)(10-1) - 89 =109. \]

The next example shows how even large examples can be done easily by hand using the recursion of Theorem \ref{Main}. 
\begin{empl} 
\label{empl:exmG}
\text{Consider $\Gamma = \langle122, 281\rangle$. We first compute the numbers $\{p_i, k_i, s \}$:}
\begin{align}
281 & = 2 \cdot 122 + 37 \nonumber \\
122 & = 3 \cdot 37 + 11 \nonumber \\
37 & = 3 \cdot 11 + 4 \nonumber \\
11 & = 2 \cdot 4 + 3 \nonumber \\
4 &= 1 \cdot 3 + 1 \nonumber 
\end{align}   
so $s = 5$, $\{p_i\} = \{122, 37, 11, 4, 3, 1 \}$ and $\{k_i\} = \{2, 3, 3, 2, 1, 3 \}$ for $0 \leq i \leq 5$. Next   for $1 \leq i \leq 5$ compute via (\ref {eqn:Ai}) the values $\{ A_i \} = \{1, 3, 10, 23, 33  \}$ and $\{ B_i \} = \{2, 7, 23, 53, 76  \}$, and by (\ref{eqn:ns}) find $\{n_i \} = \{3, 3, 0, 1, 3 \}$ and $\{N_i  \} = \{10, 7, 4, 4, 3  \}$. Next compute $\gamma_i$ and $r_i$ for $1 \leq i \leq 9$ $(= N_1-1)$, and obeying the inequalities in the theorem; from these values we immediately calculate the generators, starting with $g_1 = p + m + 1$: 
\begin{table}[H]
  \begin{center}
    \label{tab:table1}
    \begin{tabular}{c | c | r c | c | c | c  }
          $i$ & $j$ & $\gamma_i$ & &   $r_i=p_j$ & $g_i$ &  $u_i$  \\
      \hline
      1  & 5 & $(B_5-1)p \;=$ & 9150 & 1 & 404 & 9554 \\
      2  & 5 & $p \;=$ & 122 & 1  & 9555 & 9677  \\
      3  & 4 & $(A_4-1)m \;=$ & 6182 & 3 & 9678 & 15860  \\
      4  & 2 & $m \;=$ & 281 & 11 & 15863 & 16144 \\
      5  & 2 & $(A_2-1)m \;=$ & 562 & 11 & 16155 & 16717   \\
      6  & 2 & $m \;=$ & 281 & 11 & 16728 & 17009   \\
      7  & 1 & $p \;=$ & 122 & 37 & 17020 &  17142  \\
      8  & 1 & $(B_1-1)p \;=$ & 122 & 37 & 17179 & 17301   \\
      9  & 1 & $p \;=$ & 122 & 37 & 17338 &  17460  \\
    \end{tabular}
  \end{center}
\end{table} \noindent
Notice that $g_9 = g_6 + 5p$. The last generator is therefore redundant (in fact $u_9 - c_9$ exceeds $\mu$). Including $g_{-1} = p$ and $g_0 = m$, we see that $\Lambda_{\gen}$ has 10 (= $N_1$) generators, including two belonging to $\Gamma$. 

In order to relate the generators to the Tjurina number, we write them in two ways: numerically and in standard form (\ref{Delm}). To compute the standard form, apply (\ref{eqn:pi}) to each $p_j$ and use the recursion of the theorem, {\em i.e.} $g_1 = p + m + 1$ and $g_{i+1} = g_i + \gamma_i + p_j$. 
\begin{table}[H]
  \begin{center}
    \label{tab:table1}
    \begin{tabular}{c|rl||c|rl|} 
      $i$ & $pm - \alpha_im - \beta_ip =$ & $g_i$    & $i$ & $pm - \alpha_im - \beta_ip =$ & $g_i$     \\
      \hline
      1  & $pm - 88m - 75p = $ &  404  & 5  & $pm - 25m - 91p = $ &  16155 \\
      2  & $pm - 55m - 76p = $ &  9555 &6 & $pm - 26m - 84p = $ &  16728  \\
       3  & $pm - 22m - 151p = $ &  9678  &7 & $pm - 28m - 77p = $ &  17020\\
      4  & $pm - 23m - 98p = $ &  15863  &8 & $pm - 27m - 78p = $ &  17179\\
    \end{tabular}
  \end{center}
\end{table} \noindent
We can now easily compute the Tjurina number for $\Lambda_{\gen}$. Consider the set of integer pairs $(\alpha_i, \beta_i)$ induced by the generators $g_i$ in standard form, $1 \leq i \leq 8$. We re-index the pairs so that $\alpha_i < \alpha_{i+1}$. Then by necessity the second coordinate is decreasing:
\[
\{(\alpha_i, \beta_i) \} = \{(22, 151), (23, 98), (25, 91), (26, 84), (27, 78), (28, 77), (55, 76), (88, 75) \}
\]
These first-quadrant points enclose rectangles of strictly positive values $(\alpha, \beta) \leq (\alpha_i, \beta_i)$ which represent elements of $\Lambda_{gen} - \Gamma^*$, and all such elements are represented in this way. Thus the sum of areas of these rectangles counts $|\Lambda_{gen} - \Gamma^*| = \mu - \tau_{\gen}$. Letting $(\alpha_0, \beta_0) = (0, 0)$, the area sum is $\sum_1^8 (\alpha_i - \alpha_{i-1})\beta_i$. Thus
\[
 \mu - \tau_{\gen} = 22 \cdot 151 + 1 \cdot 98 + 2 \cdot 91 + 1 \cdot 84 + 1 \cdot 78 + 1 \cdot 77 + 27 \cdot 76 + 33 \cdot 75 = 8368.
\]
 \end{empl}

\section{A formula for the minimal Tjurina number} \label{Tjurinaformula}

In \cite{ABM} an explicit formula is given for the minimal (i.e. \emph{generic}) Tjurina number $\tau_{\gen}$ of any irreducible plane curve germ. We present the output of this formula in the case of a 2-generator semigroup, in terms of the notions defined in this article.

\begin{theorem} \label{eqn:Tjur} 
Let $\Gamma = \langle p,m \rangle$ as above, $\mu = (p-1)(m-1)$ the conductor of $\Gamma$, and denote by $\lfloor x \rfloor$ the floor of $x$. Then 
\begin{equation}
\tau_{\gen} = \mu - \Bigl \lfloor \frac{m}{p} \Bigr\rfloor \Bigl\lfloor \frac{(p-1)^2}{4} \Bigr\rfloor + \Bigl\lfloor \frac{p-1}{2} \Bigr\rfloor + \Bigl\lfloor \frac{p_1}{2} \Bigr\rfloor - \sum_{i=1}^{s-1} \Bigl\lfloor \frac{p_{i-1}}{p_i} \Bigr\rfloor \Bigl\lfloor \frac{p_i^2}{4} \Bigr\rfloor.
\end{equation}
\end{theorem}

The theorem follows easily from the more general formula of \cite{ABM} and the fact that the multiplicity sequence of $\Gamma$ is read off from the Euclidean algorithm. Note that the coefficients in the sum are the same as the numbers $k_i$, i.e.
\[ \Bigl\lfloor \frac{p_{i-1}}{p_i} \Bigr\rfloor = k_i.\]

\begin{empl} 
In Example~\ref{empl:exmG} we computed for the generic value set of
semigroup $\Gamma = \langle122, 281\rangle$
the formula
$
 \mu - \tau_{\gen} = 8368
$.
On the other hand, Theorem~\ref{eqn:Tjur}  calculates 
\begin{align*}
\mu - \tau_{\gen} &= \Bigl \lfloor \frac{m}{p} \Bigr\rfloor \Bigl\lfloor \frac{(p-1)^2}{4} \Bigr\rfloor - \Bigl\lfloor \frac{p-1}{2} \Bigr\rfloor - \Bigl\lfloor \frac{p_1}{2} \Bigr\rfloor + \sum_{i=1}^4 \Bigl\lfloor \frac{p_{i-1}}{p_i} \Bigr\rfloor \Bigl\lfloor \frac{p_i^2}{4} \Bigr\rfloor \\
	 &= \Bigl \lfloor \frac{281}{122} \Bigr\rfloor \Bigl\lfloor \frac{121^2}{4} \Bigr\rfloor - \Bigl\lfloor \frac{121}{2} \Bigr\rfloor - \Bigl\lfloor \frac{37}{2} \Bigr\rfloor + \Bigl\lfloor \frac{122}{37} \Bigr\rfloor \Bigl\lfloor \frac{37^2}{4} \Bigr\rfloor \\  
	 &\phantom{=-}+ \Bigl\lfloor \frac{37}{11} \Bigr\rfloor \Bigl\lfloor \frac{11^2}{4} \Bigr\rfloor + \Bigl\lfloor \frac{11}{4} \Bigr\rfloor \Bigl\lfloor \frac{4^2}{4} \Bigr\rfloor + \Bigl\lfloor \frac{4}{3} \Bigr\rfloor \Bigl\lfloor \frac{3^2}{4} \Bigr\rfloor \\
	 &= 2 \cdot 3660 - 60 - 18 + 3 \cdot 342 + 3 \cdot 30 + 2 \cdot4 + 1 \cdot 2 = 8368.
\end{align*}
\noindent
We have not yet investigated the path that connects these two approaches to calculating the generic Tjurina number. 
\end{empl} 

The authors wish to thank Patricio Almir\'on, Gary Kennedy, and Richard Montgomery for several useful conversations regarding this work.

\end{document}